\documentclass{article}
\usepackage{graphicx}
\usepackage{amsmath}
\usepackage{amssymb}

\begin{document}
\title{Explicit Bounds for the Distribution Function of the Sum of Dependent Normally Distributed Random Variables}
\author{Walter Schneider\thanks{email: schneider@fos-bos-passau.de}}
\maketitle
\pagestyle{headings}

   
\begin{abstract}
In this paper an analytic expression is given for the bounds of the distribution function of the sum
of dependent normally distributed random variables. Using the theory of copulas and the important Fr\'{e}chet bounds the dependence structure is not restricted to any
specific type. Numerical illustrations are provided to assess the
quality of the derived bounds.

\end{abstract}

\section{Introduction}
\label{intro}
Many problems in mathematical probability theory involve the computation of the distribution function for the sum of random variables (RVs). In case of independent RVs $X$ and $Y$ the distribution $F_Z$ of the sum $Z=X+Y$ is given by the integral
\begin{equation} \label{eq1}
F_Z(t)=P(X+Y \le t)=\int_{-\infty}^t \left(\int f(z-y)g(y)dy \right)dz
\end{equation}
with $f$ and $g$ as the corresponding density functions of $X$ and $Y$ respectively.

If the RVs however are stochastically dependent, which is a common situation in practice, the convolution of the marginal densities $f$ and $g$ in the integrand of \eqref{eq1} is no longer valid. Under the assumption that the joint density function $f_{XY}$ of $X$ and $Y$ is specified, the distribution $F_Z$ can be calculated by
\begin{equation} \label{eq2}
F_Z(t)=P(X+Y \le t)=\int_{-\infty}^t \left(\int f_{XY}(z-y,y)dy \right)dz
\end{equation}

An analytical exact solution of the integral both in \eqref{eq1} and \eqref{eq2} is feasible for certain marginals $f, g$ or joint densities $f_{XY}$ (e.g. the normal one). In situations where $f, g$ or $f_{XY}$ is too cumbersome to work with, one could use Monte-Carlo simulation for numerical solution of the integrals.

In many circumstances, however, the joint density of the RVs $X$ and $Y$ is not known, despite of given marginal distributions.
This leads to the question, whether it is possible to provide bounds for the distribution of the sum, which are valid for all possible joint distributions.
The original problem was formulated by A. N. Kolmogorov: Let $X$ and $Y$ be RVs with given distributions $F_X$ and $F_Y$. Find bounds $G^\wedge$ (upper bound) and $G^\vee$ (lower bound) for the distribution $G$ of the sum $Z=X+Y$, such that
\begin{eqnarray} \label{eq3}
G^\vee(z)=\inf P(X+Y<z) \\
G^\wedge(z)=\sup P(X+Y<z)
\end{eqnarray}
where the infimum and supremum are taken over all possible joint distributions having the marginal distributions $F_X, F_Y$.
In this situation, it is said, that the joint distribution has fixed margins.

By now a rich literature is available on this subject. G.D. Makarov \cite{RefA} solved Kolmogorov's problem via a cumbersome, ad hoc argument. Other authors, e.g. Frank et al. \cite{RefB} have applied theory orginally studied by Fr\'{e}chet which leads to copulas naturally. The authors in \cite{RefE}, \cite{RefF} provide examples of distributions for which the bounds $G^\vee$ and $G^\wedge$ can be explicitly computed. These distributions include the uniform, the Cauchy and the exponential families. For further reading the interested reader is referred to \cite{RefG} - \cite{RefK}.

In this contribution we compute explicit bounds for normally distributed RVs by means of the copula based theory.
Compared to the work of Frank et al. \cite{RefB} who addressed this problem 1987 by using the method of Langrange multipliers we will use a simple variable substitution. 

The rest of the paper is organized as follows. In Section 2 we give a brief outline of the copula theory and present the general theorem for bounding the distribution of the sum of RVs with unknown dependence. Then in Section 3 lower and upper bounds for the distribution of the sum of normally distributed RVs are presented. Finally we illustrate these bounds by numerical examples in Section 4.

\section{Preliminaries}
\subsection{\bf Copulas}
Copulas are used in probability theory and statistics for modeling the dependence between RVs $X_1,X_2,...,X_d$ at a deeper level allowing to understand dependence measures different from a simple correlation coefficient approach. For any natural $d$, a ($d$-dimensional) copula is a distribution function on $[0,1]^d$ with standard uniform margins.

For the formal definition of a copula as well as a general introduction into the wide fields of copulas the reader is referred to \cite{RefC}, \cite{RefD}.
In this paper the following theorems are of particular importance. \vspace{0.1cm}
\newtheorem{theoreme}{Theorem}
\begin{theoreme}[Sklar's Theorem]
Let $H$ be a joint distribution function for the RVs $X$ and $Y$ with marginals $F$ and $G$. Then there exists a copula $C$ such that for all $x, y \in \mathbb{R}$:
\begin{equation}\label{eq4}
H(x,y)=C(F(x),G(y))
\end{equation}
\end{theoreme}

\begin{theoreme}[Fr\'{e}chet Bounds]
Let $(u,v) \in I^2=([0,1],[0,1]) \subset \mathbb{R}^2$, then for every copula $C$ and every $(u,v) \in I^2$ the Fr\'{e}chet-bounds inequality is valid:
\begin{equation} \label{eq5}
W(u,v) \le C(u,v) \le M(u,v)
\end{equation}
with $M(u,v) = min(u,v)$ and $W(u,v) = max(u+v-1,0)$ as the Fr\'{e}chet bounds.
\end{theoreme}
\vspace{0.2cm}
The Fr\'{e}chet bounds inequality together with Sklar's Theorem leads to the Fr\'{e}chet bounds for the joint distribution function $H$, 
\begin{equation}\label{eq6}
max(F(x)+G(y)-1,0) \le H(x,y) \le min(F(x),G(y))
\end{equation}
so the joint distribution $H$ is bounded in terms of its own marginals.
The proofs of these theorems as well as \eqref{eq6} can be found in \cite{RefC}.

\subsection{\bf Bounding the sum of RVs with unknown dependence}
The problem of calculating the distribution $G$ of the sum $Z=X+Y$, where the dependence{\footnote{\scriptsize A note on terminology: The term "dependence" will be used for a common measure in the study of the dependence betweeen RVs and therefore is not restricted to a measure of the linear dependence between RVs.}} between $X$ and $Y$ is unknown, is called the Kolmogorov problem. Using the theory of copulas the following theorem and its proof are given in \cite{RefC}.
\begin{theoreme}
Let $X$ and $Y$ be RVs with distribution functions $F_X$ and $F_Y$. Let $G$ denote the distribution function of $X+Y$. Then
\begin{equation} \label{eq7}
G^\vee(z) \le G(z) \le G^\wedge(z)
\end{equation}
where
\begin{equation} \label{eq8}
G^\vee(z)=\sup_{x+y=z}\{W(F_X(x),F_Y(y))\}
\end{equation}
\begin{equation} \label{eq9}
G^\wedge(z)=\inf_{x+y=z}\{\tilde{W}(F_X(x),F_Y(y))\}
\end{equation}
with $\tilde{W}(u,v)=u+v-W(u,v)=min(u+v,1)$ and $W(u,v)$ as in Theorem 2.
\end{theoreme}

In the following $G^\wedge$ is named as upper bound and $G^\vee$ as lower bound for the distribution $G$.

\section{Bounds for the sum of normally distributed RVs}
Let $X$ and $Y$ be normally distributed with means $\mu_X, \mu_Y$ and standard deviations $\sigma_X, \sigma_Y$, denoted by $X \sim N(\mu_X,\sigma_X^2), Y \sim N(\mu_Y,\sigma_Y^2)$. Their distribution functions $F_X$ and $F_Y$ are given by
\begin{displaymath}
F_X(x)=\Phi\left(\frac{x-\mu_X}{\sigma_X} \right) \hspace{0.5cm} \textnormal{and} \hspace{0.5cm} F_Y(y)=\Phi\left(\frac{y-\mu_Y}{\sigma_Y} \right)
\end{displaymath}
where $\Phi$ denotes the standard normal distribution
\begin{displaymath}
\Phi(t)=\frac{1}{\sqrt{2\pi}} \int_{-\infty}^t e^{-\frac{s^2}{2}}ds
\end{displaymath}

\vspace{0.3cm}
In case of $\sigma_X \neq \sigma_Y$ the lower and upper bound for the distribution $G$ of the sum $Z=X+Y$ are given by the following proposition. 

\newtheorem{proposition}{Proposition}
\begin{proposition}
The lower bound $G^\vee$ and the upper bound $G^\wedge$ for $G$ are calculated by
\begin{equation} \label{eq10}
G^\vee(z)=\left\{\begin{array}{ll} \max \left(\Psi(x_1)-1,0\right): & \hspace{0.5cm} \Psi(x_1) > \Psi(x_2) \\
\max \left(\Psi(x_2)-1,0 \right): & \hspace{0.5cm} \Psi(x_1) < \Psi(x_2) \end{array} \right.
\end{equation}
and
\begin{equation} \label{eq11}
G^\wedge(z)=\left\{\begin{array}{ll} \min \left(\Psi(x_1),1 \right): & \hspace{0.5cm} \Psi(x_1) < \Psi(x_2) \\
\min \left(\Psi(x_2),1 \right): & \hspace{0.5cm} \Psi(x_1) > \Psi(x_2) \end{array} \right.
\end{equation}
with the function $\Psi:\mathbb{R} \rightarrow \mathbb{R}$:
\begin{equation} \label{eq12}
\Psi(x)=\Phi\left(\frac{x-\mu_X}{\sigma_X}\right)+\Phi\left(\frac{z-x-\mu_Y}{\sigma_Y}\right)
\end{equation}
and $x_1, x_2$ as the local extremas of $\Psi$.
\end{proposition}
{\it{Proof.}} 
Using Theorem 3 together with the Fr\'{e}chet bounds and \eqref{eq6} introduced in section 2, $G^\vee$ and $G^\wedge$ are given by
\begin{eqnarray} 
G^\vee(z)&=&\sup_{x+y=z} \max\left\{\Phi\left(\frac{x-\mu_X}{\sigma_X}\right)+\Phi\left(\frac{y-\mu_Y}{\sigma_Y}\right)-1,0\right\} \label{eq13} \\
G^\wedge(z)&=&\inf_{x+y=z} \min\left\{\Phi\left(\frac{x-\mu_X}{\sigma_X}\right)+\Phi\left(\frac{y-\mu_Y}{\sigma_Y}\right),1\right\} \label{eq14}
\end{eqnarray}
For the sum of the distributions $F_X$ and $F_Y$ in \eqref{eq13}, \eqref{eq14} we introduce the function $\Psi:\mathbb{R} \rightarrow \mathbb{R}$:
\begin{equation} \label{eq15}
\Psi(x)=\Phi\left(\frac{x-\mu_X}{\sigma_X}\right)+\Phi\left(\frac{z-x-\mu_Y}{\sigma_Y}\right)
\end{equation}
where the variable $y$ has been substituted by $z-x$.

Now, we can proceed with the computation of the local extremas of $\Psi$ as function of one variable. Therefore the first derivative $\Psi'(x)$ is built:
\begin{equation} \label{eq16}
\Psi'(x)=\frac{1}{\sigma_X\sqrt{2\pi}} \cdot e^{-\frac{(x-\mu_X)^2}{2\sigma_X^2}}-\frac{1}{\sigma_Y\sqrt{2\pi}} \cdot e^{-\frac{(x-(z-\mu_Y))^2}{2\sigma_Y^2}}
\end{equation}
Finding the zeros of $\Psi'(x)$ leads to 
\begin{equation} \label{eq17}
\frac{1}{\sigma_X\sqrt{2\pi}} \cdot e^{-\frac{(x-\mu_X)^2}{2\sigma_X^2}}=\frac{1}{\sigma_Y\sqrt{2\pi}} \cdot e^{-\frac{(x-(z-\mu_Y))^2}{2\sigma_Y^2}} 
\end{equation}
which - after some technical calculation - is equivalent to the quadratic equation
\begin{equation} \label{eq18}
\alpha x^2 + \beta x + \gamma = 0  
\end{equation}
with the variables $\alpha, \beta, \gamma$ as follows:
\begin{eqnarray} \label{eq19}
\alpha&=&\frac{1}{2\sigma_X^2}-\frac{1}{2\sigma_Y^2} \\ \label{eq19_1}
\beta&=&\frac{z-\mu_Y}{\sigma_Y^2}-\frac{\mu_X}{\sigma_X^2} \\ \label{eq19_2}
\gamma&=&\frac{-(z-\mu_Y)^2}{2\sigma_Y^2}+\frac{\mu_X^2}{2\sigma_X^2}-\ln{\frac{\sigma_Y}{\sigma_X}} \label{eq19_3}
\end{eqnarray}
The solutions of \eqref{eq18} are given by $x_{1,2}$
\begin{equation} \label{eq20}
x_{1,2}=\frac{-\beta\pm \sqrt{\beta^2-4\alpha \gamma}}{2\alpha}
\end{equation}
which are the candidates for the local extremas of the function $\Psi(x)$. 

As $\sigma_X \neq \sigma_Y$ ($\alpha \neq 0$) it is guaranteed that the division in \eqref{eq20} is defined. Moreover, any possible values for $\mu_X, \mu_Y, \sigma_X, \sigma_Y$ will lead to two solutions $x_1 \neq x_2$ of the quadratic equation \eqref{eq18}.
In order to show this, we compute the limits of $\Psi$ for $x \rightarrow \pm \infty$.
\begin{eqnarray} 
\lim_{x \rightarrow -\infty} \Psi(x) = \lim_{x \rightarrow -\infty} \left(\underbrace{ \Phi\left(\frac{x-\mu_X}{\sigma_X}\right)}_{\rightarrow 0}+\underbrace{\Phi\left(\frac{z-x-\mu_Y}{\sigma_Y}\right)}_{\rightarrow 1} \right) = 1  \\ \label{eq21} 
\lim_{x \rightarrow \infty} \Psi(x) = \lim_{x \rightarrow \infty} \left(\underbrace{ \Phi\left(\frac{x-\mu_X}{\sigma_X}\right)}_{\rightarrow 1}+\underbrace{\Phi\left(\frac{z-x-\mu_Y}{\sigma_Y}\right)}_{\rightarrow 0} \right) = 1  \label{eq22} 
\end{eqnarray}

The same limits of $\Psi$ for $x \rightarrow \pm \infty$ together with the fact that $\Psi$ is continuous but not constant on $\mathbb{R}$ implies that $\Psi$ has one extremum at least. If $\Psi$ had only one extremum there would be exactly one solution $x_1=x_2$ (zero of order 2) of \eqref{eq18}. A zero of order 2 for $\Psi'(x)=0$ however cannot be a extreme value. So, in any case there exist exactly two extreme values.

Now we consider the set $U:=\{\Psi(x_1)-1,\Psi(x_2)-1,0\}$, which contains both a minimum and maximum value. If $\Psi(x_1) > \Psi(x_2)$ then $\max (U)=\max\left(\Psi(x_1),0\right)$, otherwise $\max (U)=\max\left(\Psi(x_2),0\right)$.  As for each set $M$ having a maximum value there is $\sup M = \max M$, equation \eqref{eq10} of Proposition 1 is shown.
Similarly equation \eqref{eq11} of Proposition 1 can be proven using the set $V:=\{\Psi(x_1),\Psi(x_2),1\}$. If $\Psi(x_1) < \Psi(x_2)$ then $\min (U)=\min\left(\Psi(x_1),1\right)$, otherwise $\min (U)=\min\left(\Psi(x_2),1\right)$. \hfill $\square$

\vspace{0.5cm}
In the special case, that the RVs have the same standard deviation ($\sigma_X=\sigma_Y$) the bounds $G^\vee$ and $G^\wedge$ can be expressed in closed form by the following corollary.
\newtheorem{corollary}{Corollary}
\begin{corollary}
Let $\sigma_X=\sigma_Y=\sigma$ the bounds $G^\vee$ and $G^\wedge$ for the distribution $G$ of the sum of two normally distributed random variables are given by
\begin{equation} \label{eq23}
G^\vee(z)=\left\{\begin{array}{ll} 0, & \hspace{0.5cm} z \leq \mu_X+\mu_Y \\ 2\Phi\left(\frac{z-\mu_X-\mu_Y}{2\sigma}\right)-1, & \hspace{0.5cm} z \geq \mu_X+\mu_Y \end{array} \right.
\end{equation}
and
\begin{equation} \label{eq24}
G^\wedge(z)=\left\{\begin{array}{ll} 2\Phi\left(\frac{z-\mu_X-\mu_Y}{2\sigma}\right), & \hspace{0.5cm} z \leq \mu_X+\mu_Y \\ 1, & \hspace{0.5cm} z \geq \mu_X+\mu_Y \end{array} \right.
\end{equation}
\end{corollary}
{\it{Proof.}}
By using the same arguments as in the proof of Proposition 1, we get $\alpha=0$ from \eqref{eq19}.
This means we have to solve a linear equation $\beta x + \gamma = 0$. From \eqref{eq19_1} and \eqref{eq19_2} we get
\begin{eqnarray}
\beta &=& \frac{z-\mu_Y}{\sigma^2}-\frac{\mu_X}{\sigma^2}=\frac{z-\mu_Y-\mu_X}{\sigma^2} \\
\gamma &=& \frac{-(z-\mu_Y)^2}{2\sigma^2}+\frac{\mu_X^2}{2\sigma^2}-\ln{1}=\frac{-z^2+2z\mu_Y-\mu_Y^2+\mu_X^2}{2\sigma^2}
\end{eqnarray}
Then the solution $x_0$ for $\beta x + \gamma = 0$ is given by
\begin{equation}
x_0=-\frac{\gamma}{\beta}=\frac{z^2-2z\mu_Y+\mu_Y^2-\mu_X^2}{2(z-\mu_Y-\mu_X)}=\frac{z-\mu_Y+\mu_X}{2}
\end{equation}
As $x_0$ is a zero of order 1 of $\Psi'(x)=0$ it is guaranteed, that $x_0$ is indeed a extreme value.
Inserting $x_0$ into $\Psi(x)$ leads to
\begin{eqnarray} \label{eq25}
\Psi(x_0)&=&\Phi\left(\frac{\frac{z-\mu_Y+\mu_X}{2}-\mu_X}{\sigma}\right)+\Phi\left(\frac{z-\frac{z-\mu_Y+\mu_X}{2}-\mu_Y}{\sigma}\right) \nonumber \\
&=&\Phi\left(\frac{z-\mu_Y-\mu_X}{2\sigma}\right)+\Phi\left(\frac{z-\mu_X-\mu_Y}{2\sigma}\right) \nonumber \\ &=&2\Phi\left(\frac{z-\mu_Y-\mu_X}{2\sigma}\right)
\end{eqnarray}
If $z>\mu_X+\mu_Y$ then the argument of the distribution function $\Phi$ in \eqref{eq25} is positive. As $\Phi(t)>0.5$ for any positive value $t \in \mathbb{R}$ we have $\Psi(x_0)>1$. Conversely, if $z<\mu_X+\mu_Y$ then the argument of the distribution function $\Phi$ in \eqref{eq25} is negative and as a result $\Psi(x_0)<1$. Together with the equations for $G^\vee$ and $G^\wedge$ in \eqref{eq13}, \eqref{eq14} the corollary is proven. \hfill $\square$

\newtheorem{remark}{Remark}
\begin{remark}
The result of corollary 1 is in accordance with the result of M. J. Frank \cite{RefB}. However the bounds from Frank in the common case for different standard deviations ($\sigma_X \neq \sigma_Y$) cannot be reproduced. Simple numerical examples have revealed, that Frank's bounds have also negative values, which is in contradiction to a distribution function.
\end{remark}

\section{Illustrative Examples}
In this section we want to illustrate the Proposition 1 in section 3 by a simple example with the following prerequisites:
\begin{itemize}
\item[$\bullet$] The RVs $X$ and $Y$ are normally distributed, where the parameters $\mu$ and $\sigma$ are chosen as follows: $X \sim N(1,0.1)$, $Y \sim N(1.5,0.15)$.  

\item[$\bullet$] For modeling the dependence amongst $X$ and $Y$ two different strategies are used: 1) $X$, $Y$ are bivariate normal distributed with a given correlation coefficient $\rho$. 2) The dependence of $X$ and $Y$ is modeled using either the Clayton ($C_{\theta}^{Cl}$) or the Gumbel ($C_{\theta}^{Gu}$) Copula \cite{RefD}. \\
\begin{eqnarray}
C_{\theta}^{Cl}(u_1,u_2)&=&(u_1^{-\theta}+u_2^{-\theta}-1)^{-\frac{1}{\theta}}  \nonumber \\
C_{\theta}^{Gu}(u_1,u_2)&=&exp{\left(-((-ln (u_1))^\theta + (-ln (u_2))^\theta)^\frac{1}{\theta} \right)} \nonumber
\end{eqnarray}
These copulas have been chosen, because the Gumbel (Clayton) copula turns out to have upper (lower) tail dependence. 
A method for generating realizations from a particular copula within Matlab is described in \cite{RefM}. The value $\theta$ has been set to $2.5$ for both copulas.

\item[$\bullet$] The bound computation $G^\vee$ and $G^\wedge$ as well as generating dependent random variables using the concept of copulas is done via Matlab.
\end{itemize}

\begin{figure}[h]
\includegraphics[width=13.0cm]{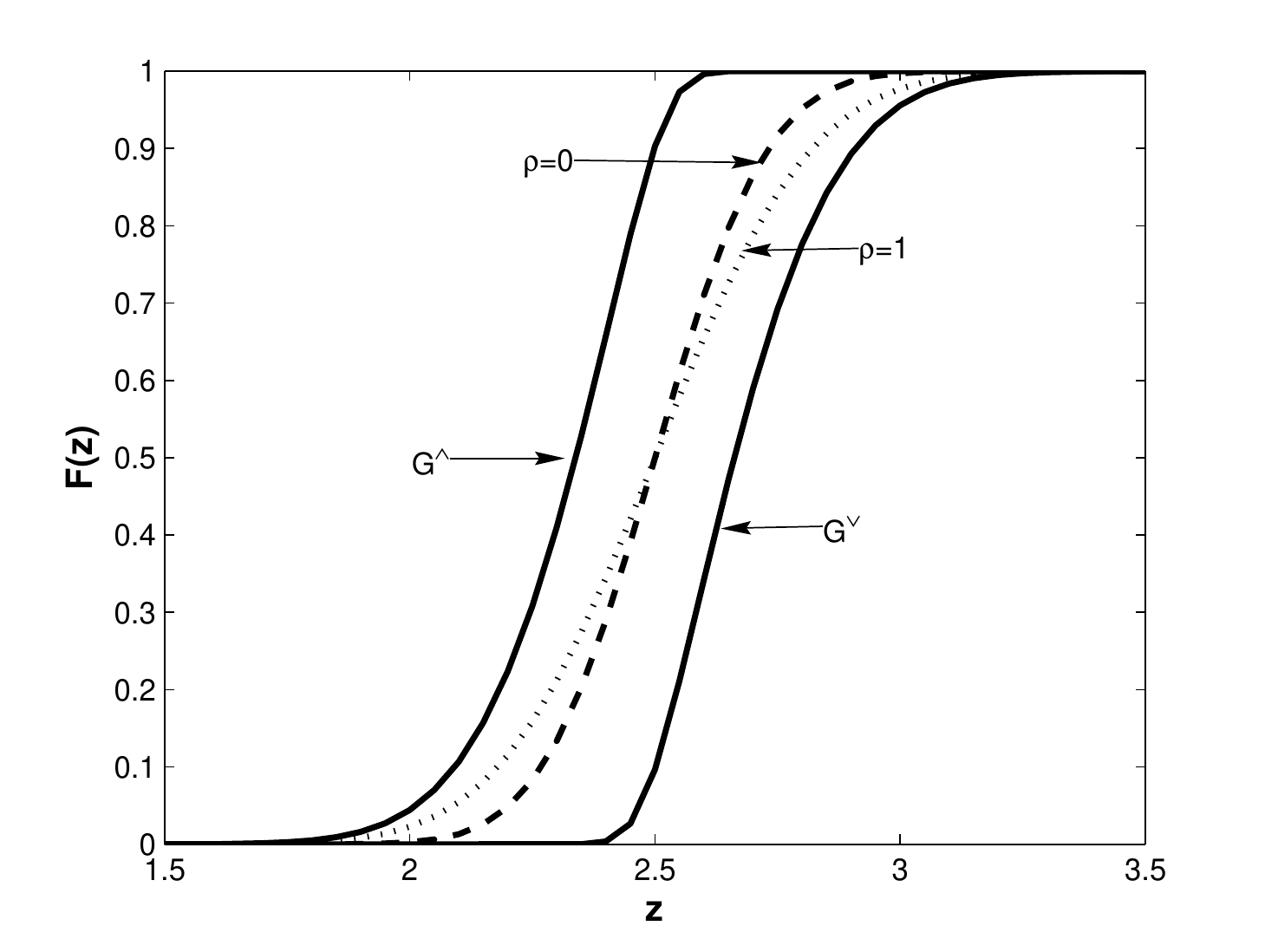}
\caption{Computed bounds $G^\wedge$, $G^\vee$ compared to the distribution $F$ of the sum $X+Y$ assuming bivariate normally distributed random variables $X, Y$}
\label{fig:eins}       
\end{figure}

In figure~\ref{fig:eins} the bounds $G^\vee$, $G^\wedge$ are compared to the sum distribution $F$ assuming $X, Y$ are bivariate normally distributed. The corner cases $\rho=0$ and $\rho=1$ have been used for the bivariate normal distribution.

\begin{figure}[h]
\includegraphics[width=13.0cm]{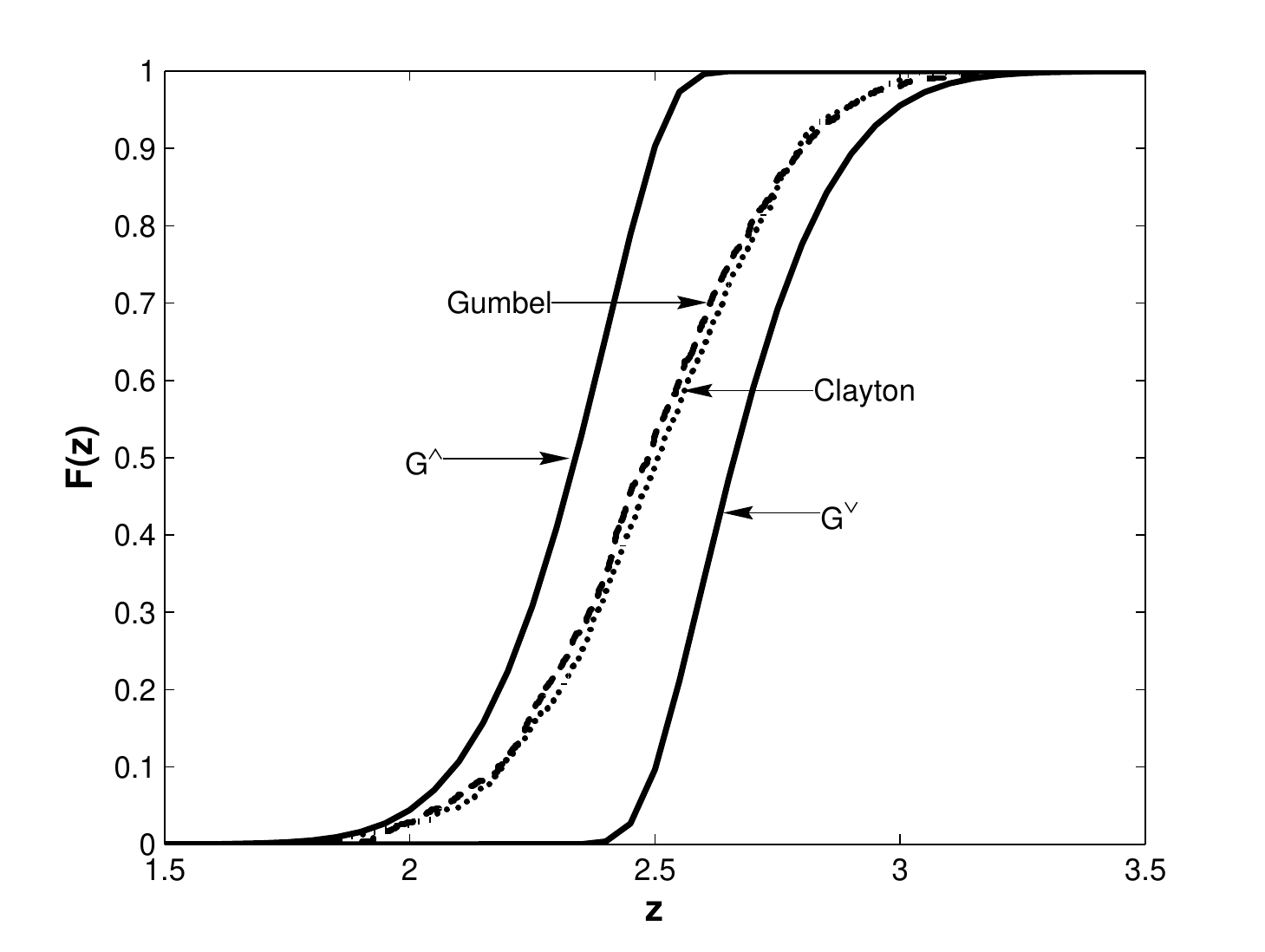}
\caption{Computed bounds $G^\wedge$, $G^\vee$ compared to the distribution $F$ of the sum $X+Y$ where the dependence of $X$ and $Y$ is modeled by Clayton and Gumbel copula}
\label{fig:zwei}       
\end{figure}

In figure~\ref{fig:zwei} the bounds $G^\vee$, $G^\wedge$ are compared to the sum distribution $F$ assuming that $X, Y$ are dependent either by Clayton or Gumbel copula.

\end{document}